\documentclass[12pt,a4paper]{article}

\usepackage{amsthm,amsmath,natbib,amssymb,amsfonts,bm}
\usepackage{epsfig}
\usepackage{placeins}
\usepackage{epsf}
\RequirePackage[dvips]{hyperref}

\newtheorem*{theorem*}{Theorem}

\newcommand{\bX}{\boldsymbol{X}}

\theoremstyle{definition}

\theoremstyle{remark}

\begin{document}

\baselineskip=21pt

\begin{center}
 {\bf \Large On randomized confidence intervals for the binomial probability}
\end{center}

\bigskip

\begin{center}
{\bf \large Paul Kabaila$^{\textstyle{^*}}$}
\end{center}

\medskip

\noindent{\it $\textstyle{^*}$Department of Mathematics and
Statistics, La Trobe University, Victoria 3086, \newline Australia}

\bigskip
\noindent{\bf Abstract}
\medskip

Suppose that  $X_1,X_2,\ldots,X_n$ are independent and
identically Bernoulli($\theta$) distributed.
Also suppose that our aim is to find an exact confidence interval for $\theta$
that is the intersection of a $1-\alpha/2$ upper confidence interval and a
$1-\alpha/2$ lower confidence interval. The Clopper-Pearson interval is the
standard such confidence interval for $\theta$, which is widely used in practice.
We consider the randomized confidence interval
of Stevens, 1950 and present some extensions, including pseudorandomized
confidence intervals. We also consider the
``data-randomized'' confidence interval of Korn, 1987 and point
out some additional attractive features of this interval. We also
contribute to the discussion about the practical use of such
confidence intervals.

\bigskip
\bigskip


\vbox{\vskip 4cm}


\noindent $^*$ Corresponding author. Address: Department of
Mathematics and Statistics, La Trobe University, Victoria 3086,
Australia; Tel.: +61-3-9479-2594; fax: +61-3-9479-2466. \newline
{\it E-mail address:} P.Kabaila@latrobe.edu.au.

\newpage

\noindent {\bf 1. Introduction}

\medskip

Suppose that  $X_1, X_2,\ldots,X_n$ are independent and
identically distributed (iid), each with a Bernoulli($\theta$) distribution ($\theta \in [0,1]$).
Let $\bX = (X_1, \dots, X_n)$. Our objective is to find a confidence interval
for $\theta$ of the form $\big[ \ell(\bX, V), u(\bX, V) \big]$, where the interval endpoints may
depend on an auxiliary random variable $V$, such that
both of the following conditions are satisfied:
\begin{align}
\label{upper_tail}
P_{\theta} \big(\theta > u(\bX, V) \big ) &\le \frac{\alpha}{2} \ \ \text{for all} \ \ \theta \\
\label{lower_tail}
P_{\theta} \big(\theta < \ell(\bX, V) \big ) &\le \frac{\alpha}{2} \ \ \text{for all} \ \ \theta.
\end{align}
Of course, such a confidence interval has infimum coverage probability that is greater than or
equal to $1-\alpha$.
The conditions \eqref{upper_tail} and \eqref{lower_tail} make the endpoints
of the confidence interval, $\ell(\bX, V)$ and $u(\bX, V)$, easy to interpret.
Confidence intervals that satisfy these conditions are the discrete-data analogue of an equi-tailed
confidence interval based on continuous data.
The solution favoured by statistical practitioners is
to find a non-randomized confidence interval for $\theta$
based on $Y = X_1 + X_2 + \dots + X_n$, which has a Binomial$(n, \theta)$
distribution. The resulting Clopper-Pearson interval (Clopper and Pearson, 1934)
is widely used in practice.
Of course, if the conditions \eqref{upper_tail} and \eqref{lower_tail} were to be replaced by the less
stringent requirement that
$P_{\theta} \big(\ell(\bX, V) \le \theta \le u(\bX, V) \big ) \ge 1 - \alpha$ for all $\theta$
then other non-randomized confidence intervals such as that of Blaker (2000) would come into consideration.
Nonetheless, there is still a lively interest in randomized and related confidence intervals, as evidenced by
e.g. Geyer and Meeden (2005) and the resulting published comments.
In the present paper, we will compare various randomized, ``pseudorandomized'' and ``data-randomized''
confidence intervals that satisfy \eqref{upper_tail} and \eqref{lower_tail} with the Clopper-Pearson confidence interval (described, for the
reader's convenience, in Section 2).

A randomized confidence interval can be found by
considering the artificial data $Z = Y + V$, where $V$ and $Y$ are independent
random variables and $V$ has a uniform distribution on $(0,1)$ (Stevens, 1950).
Equi-tailed $1-\alpha$ confidence intervals based on $Z$ dominate the
$1-\alpha$ Clopper-Pearson confidence intervals. In Section 3, we review this randomized confidence interval and introduce
some extensions to the idea of a randomized confidence interval, including pseudorandomized
confidence intervals. In Section 4, we review the usual
objections to the use of randomized confidence intervals in practice and note a new
objection based on the need to condition on an ancillary statistic.

Korn (1987) introduced a ``data-randomized'' confidence interval for $\theta$ that
uses the data itself to generate the randomization. This confidence interval
does not require the
use of an auxiliary variable $V$ and overcomes some of the objections to the
use of randomized confidence intervals in practice. In Section 5, we
review this confidence interval and point out some
additional attractive features of this interval. In Section 6, we note the objections
of Senn (2007ab) to the use of such intervals in practice and note a further objection based
on an invariance argument.

In Section 7, we consider the properties of an unusual confidence interval for $\theta$
that turns out to be a ``data-randomized'' confidence interval. We also explain why
we expect this confidence interval to have properties that are inferior to the
``data-randomized'' confidence interval of Korn (1987).

\bigskip

\noindent {\bf 2. Clopper-Pearson confidence interval for the binomial
probability}

\medskip

Let $f_{\theta}(y) = P_{\theta}(Y=y)$ and $F_{\theta}(y) = P_{\theta}(Y \le y)$.
For observed value $y$ of $Y$, the Clopper-Pearson $1-\alpha$ confidence interval for
$\theta$ is found as follows.
The p-value for testing the null hypothesis $H_0: \theta = \tilde{\theta}$ against the
alternative hypothesis $H_A: \theta < \tilde{\theta}$ is $P_{\tilde{\theta}}(Y \le y)$.
A $1-\alpha$ upper confidence interval for $\theta$ is $\big \{\theta: P_{\theta}(Y \le y) > \alpha/2 \big \}$.
We replace all upper confidence intervals of the form $[b,c)$ by $[b,c]$.
This does not decrease the coverage probability of this confidence interval and leaves its length
unchanged.
The upper endpoint of the Clopper-Pearson $1-\alpha$ confidence interval is 1 for $y=n$; otherwise it is
the solution for $\theta$ of
\begin{equation}
\label{upper_endpoint_eqn_CP}
P_{\theta}(Y \le y) = F_{\theta}(y) = \frac{\alpha}{2}.
\end{equation}
The lower endpoint of this interval is 0 for $y=0$; otherwise it is the solution
for $\theta$ of
\begin{equation}
\label{lower_endpoint_eqn_CP}
P_{\theta}(Y \ge y) = 1 - F_{\theta}(y-1) = \frac{\alpha}{2}.
\end{equation}
Convenient expressions for the solutions of these equations
are described e.g. by Casella and Berger, (2002, p.454). We denote the
Clopper-Pearson interval by $\big[\ell_{CP}(Y), u_{CP}(Y) \big]$.



\bigskip

\noindent {\bf 3. Randomized confidence interval for the binomial
probability}

\medskip

The randomized confidence interval of Stevens (1950) can be found by considering
the artificial data $Z = Y + V$, where $V$ and $Y$ are independent random variables
and $V \sim U(0,1)$.
Also assume that either $V \in [0,1)$ or $V \in (0,1]$.
For observed values $y$, $v$ and $z$ of $Y$, $V$ and $Z$, respectively,
this confidence interval for $\theta$ is found as follows.
The p-value for testing $H_0: \theta = \tilde{\theta}$ against
$H_A: \theta < \tilde{\theta}$ is $P_{\tilde{\theta}}(Z \le z)$.
A $1-\alpha$ upper confidence interval for $\theta$ is $\big \{\theta: P_{\theta}(Z \le z) > \alpha/2 \big \}$.
We replace all upper confidence intervals of the form $[b,c)$ by $[b,c]$.
This does not decrease the coverage probability of this confidence interval and leaves its length
unchanged.
The upper endpoint of
the randomized interval is 1 for $y=n$ and $v > \alpha/2$ ; otherwise it is
the solution for $\theta$ of
\begin{equation}
\label{upper_endpoint_eqn_randomized}
P_{\theta}(Z \le z)
= v f_{\theta}(y) + F_{\theta}(y-1)
= (1-v) F_{\theta}(y-1) + v F_{\theta}(y)
= \frac{\alpha}{2}.
\end{equation}
The lower endpoint of this interval is 0 for $y=0$ and $v < 1-\alpha/2$;
otherwise it is the solution
for $\theta$ of
\begin{equation}
\label{lower_endpoint_eqn_randomized}
P_{\theta}(Z \ge z)
= (1-v) f_{\theta}(y) + 1 - F_{\theta}(y)
= (1-v) \big(1-F_{\theta}(y-1) \big) + v \big(1-F_{\theta}(y)\big)
= \frac{\alpha}{2}.
\end{equation}
Denote the resulting confidence interval by $\big[\ell_R(y,v), u_R(y,v) \big]$.
This interval satisfies both of the following conditions:
\begin{align*}
P_{\theta} \big(\theta > u_R(Y,V) \big ) &= \frac{\alpha}{2} \ \ \text{for all} \ \ \theta \\
P_{\theta} \big(\theta < \ell_R(Y,V) \big ) &= \frac{\alpha}{2} \ \ \text{for all} \ \ \theta.
\end{align*}
By comparing the left-hand sides of \eqref{upper_endpoint_eqn_randomized}
with \eqref{upper_endpoint_eqn_CP} and \eqref{lower_endpoint_eqn_randomized}
with \eqref{lower_endpoint_eqn_CP}, we find that for every $y$ and $v \in [0,1]$,
this confidence interval is contained strictly within the Clopper-Pearson $1-\alpha$
confidence interval. The confidence interval lower endpoint $\ell_R(y,v)$ is a
nondecreasing function of $v$ that is strictly increasing
 for (a) $y=0$ and $v \ge 1-\alpha/2$ and (b) all $y \ge 1$.
The confidence interval upper endpoint $u_R(y,v)$ is
nondecreasing function of $v$ that is strictly increasing
for (a) $y=n$ and $v \le \alpha/2$ and (b) all $y \le n-1$.

The confidence interval $\big[\ell_R(Y,V), u_R(Y,V) \big]$ dominates the
$1-\alpha$ Clopper-Pearson confidence interval.
This domination is
possible because ``losses for interval estimation and hypothesis testing are not
usually convex'' (Casella and Berger, 1999, p.484).
For the $1-\alpha$ Clopper-Pearson confidence interval $\big[\ell_{CP}(Y), u_{CP}(Y) \big]$,
$P_{\theta} \big(\theta > u_{CP}(Y) \big )$ and $P_{\theta} \big(\theta < \ell_{CP}(Y) \big )$
are discontinuous functions that typically take values well above $\alpha/2$ for some
values of $\theta$. By contrast, the confidence interval $\big[\ell_R(Y,V), u_R(Y,V) \big]$
has ideal coverage properties. The excellent theoretical properties of this randomized
interval can be traced to the fact that the addition of $V$ to $Y$ has ``split''
each observation $y$ into a continuous set of values, where
the values that $y$ is split into are less than all of the values that $y+1$ is split into
for each $y=0, \dots, n-1$. In the language of
Kabaila and Lloyd (2006), $Z$ is a ``refinement'' of $Y$.

As described in Appendix A, the confidence interval $\big[\ell_R(Y,V), u_R(Y,V) \big]$
may be generalized by allowing the lower and upper endpoints to depend on different
random variables $V_{\ell}$  and $V_u$, respectively, where each of these random
variables is uniformly distributed on $(0,1)$. However, as explained in Appendix A,
there seems to be no advantage to be gained from this generalization.

As described in Appendix B,
we may also construct a randomized confidence interval for $\theta$ using an auxiliary discrete random variable
$W$. This random variable may be viewed as an approximation to $V$, which has a uniform
distribution on $(0,1)$.

The usual interpretation of the coverage probability of a confidence interval is that, in a sequence of
independent repetitions of the statistical experiment that gave rise to this confidence interval,
the long-run proportion of confidence intervals that includes the parameter is equal to the coverage
probability. As described in Appendix C, this interpretation allows us to consider confidence intervals for $\theta$ that are
influenced by an appropriately-chosen auxiliary deterministic sequence, instead of the observed value of an
auxiliary random variable such as $V$. These deterministic sequences may be pseudorandom, quasi-random
or possess a very obvious pattern. What we do in this appendix is to replace expectations by the
corresponding long-run averages.


\bigskip

\noindent {\bf 4. Objections to the use of randomized confidence intervals in practice}

\medskip

Cox and Hinkley (1974, p.100)
view randomization of this type as ``a mathematical artifice'' that is ``of
no direct practical importance''. Two very cogent objections
to the use of randomized confidence intervals in practice are the following:

\smallskip

\noindent (1) Two scientists using the same procedure to construct
a randomized $1-\alpha$ confidence interval for $\theta$ based on the same
observed value $y$ will, with probability 1, produce different confidence intervals.

\smallskip

\noindent (2) The randomized interval is influenced by an auxiliary random variable
$V$ that has no relation to the problem under consideration.

\smallskip

\noindent These two reasons are presented, for example, by Kiefer (1987, p.50) and Korn (1987, p.707).
Would the first of these objections be reduced if the following procedure were adopted?
A website maintained by a reputable organisation would, upon the provision of the name
of the user and the title of a project, provide an observation $v$ of $V \sim U(0,1)$ derived
from a genuinely random source, such as electronic thermal noise. Together with this
observation, this website would provide a identification number. The user would then use
this observation $v$ to construct his/her realisation of a randomized confidence interval
and report this interval, together with $v$ and this identification number. The website
would permanently list the names of all users, projects, identification numbers and
values of $v$.

In Appendix C, we show how confidence intervals depending on an
appropriately-chosen auxiliary deterministic sequence have the desired long-run properties.
Such a sequence may be pseudorandom, quasi-random or may be a sequence with a very obvious
pattern.
However, it would seem that the alarm experienced by practitioners in response to having their confidence
interval being influenced by an auxiliary variable increases as we move from random variable to
preudorandom variable to quasi-random variable to a variable showing a very obvious pattern.

We now add a third reason for rejecting the use of randomized confidence intervals in practice.
Statisticians who believe that inference should be carried out
conditional on an appropriate ancillary statistic (see e.g. Cox and Hinkley, 1974)
would have the following objection
to the use of such a randomized confidence interval in practice. The random
variables $Y$ and $V$ can be recovered from the random variable $Z$. The statistic
$V$ has a distribution that does not depend on $\theta$ i.e. it is an ancillary
statistic. Carrying out inference conditional on $V=v$ is equivalent to carrying
out inference based solely on $Y$, leading to a non-randomized confidence interval.

\bigskip

\noindent {\bf 5. ``Data-randomized'' confidence interval for the binomial
probability}

\medskip

An apparent solution to the first of the objections described in the previous section
and a mitigation of the
second and third objections has been proposed by Korn (1987). This author
defines $W$ to be the one-sided p-value from the Wilcoxon rank-sum test for testing
the null hypothesis that the ones in the sequence $X_1, \dots, X_n$ are
randomly distributed in this sequence against the alternative hypothesis
that they come near the beginning of this sequence. The distribution of $W$, conditional
on $Y=y$, is uniform on $\big \{1/{n \choose y}, 2/{n \choose y}, \dots, {n \choose y}/{n \choose y} \big \}$,
so that it does not depend on $\theta$. This
conditional distribution
stochastically dominates the distribution of $V \sim U(0,1)$. Korn (1987) uses this to prove
that \eqref{upper_tail_discrete}, stated in Appendix B, holds true.
%
%

Korn (1987) does not describe how the lower endpoint of his randomized confidence
interval should be found. Based on the work presented in Appendix B, it is clear that
the lower endpoint of this interval should be found as follows.
Define the discrete random
variable $\tilde{W}$ by the requirement that, conditional on $Y=y$, $\tilde{W}=W-1/{n \choose y}$. Thus, conditional on $Y=y$,
$\tilde{W}$ is uniformly distributed
on $\big\{0, 1/{n \choose y}, \dots, \big({n \choose y}-1 \big)/{n \choose y} \big\}$.
Hence the data-randomized confidence interval for $\theta$ is
$\big[\ell_R(Y,\tilde{W}), u_R(Y,W) \big]$. This interval satisfies
\eqref{upper_tail_discrete}
and
\eqref{lower_tail_discrete} (stated in Appendix B).
If $n$ is not too small then \eqref{upper_tail_discrete_approx}
and
\eqref{lower_tail_discrete_approx} (stated in Appendix B) are also satisfied
and the expected length functions of the confidence intervals
$\big[\ell_R(Y,V), u_R(Y,V) \big]$ and $\big[\ell_R(Y,\tilde{W}), u_R(Y,W) \big]$
are approximately equal.
It is straightforward to show that
$\big[\ell_R(Y,\tilde{W}), u_R(Y,W) \big]$ dominates the
$1-\alpha$ Clopper-Pearson confidence interval.
This confidence interval
eliminates the first of the objections raised in Section 4
since the data determines the randomization.
Consequently, Korn (1987) calls these ``data-randomized'' confidence intervals.

The excellent theoretical properties of this data-randomized
interval can be traced to the fact that the addition of $W$ to $Y$ has ``split''
each observation $y$ into $n \choose y$ values, where
the values that $y$ is split into are less than all of the values that $y+1$ is split into
for each $y=0, \dots, n-1$. In the language of
Kabaila and Lloyd (2006), $Y+W$ is a ``refinement'' of $Y$. The upper endpoints of the data-randomized confidence intervals
are based on $Y+W$, which can take
\begin{equation*}
\sum_{y=0}^n {n \choose y} = \sum_{y=0}^n {n \choose y} 1^{n-y} 1^y = (1+1)^n = 2^n
\end{equation*}
possible values.
The excellent theoretical
properties of this data-randomized confidence interval suggest that it can be chosen
as a standard against which other data-randomized confidence intervals can be judged.

Of course, conditional on $Y=y$, there are $n \choose y$ equally-likely distinct locations
of the $y$ ones. Any one-to-one correspondence between these distinct locations and
the integers $1, 2, \dots, {n \choose y}$ can be used, in the obvious way, to generate a random variable with
the same conditional distribution as $W$. This random variable could be used as an alternative
to $W$ to construct a data-randomized confidence interval with the same coverage and expected
length properties as $\big[\ell_R(Y,\tilde{W}), u_R(Y,W) \big]$.

\bigskip

\noindent {\bf 6. Objections to the use of data-randomized confidence intervals in practice}

\medskip

Senn (2007ab) has objected to data-randomized inference procedures on two general
grounds that specialise in the present circumstance to the following:

\smallskip

\noindent (1) The ``split'' of each observation $y$ leads to quite an arbitrary ranking
of the values into which $y$ is split. Senn (2007a) says that any such split should
be based only on some meaningful comparison of the values that arise from a given observation
$y$. In the present circumstance, there does not appear to be any meaningful
comparison that could be used as the basis for this split.

\smallskip

\noindent (2) The confidence interval described by Korn (1987) is only one of many possible
data-randomized confidence intervals with the same theoretical properties. If $W$ is
replaced by $W^*$, which is obtained by calculating $W$ after the observations
have undergone a given permutation then the resulting data-randomized confidence
intervals have the same theoretical properties. Thus users of data-randomized
confidence intervals will only be able to find a unique $1-\alpha$
confidence interval for given data $x_1, \dots, x_n$ if a {\sl convention} can be established
that the interval is based only on the auxiliary random variable $W$ proposed by Korn (1987)
(and not some alternative auxiliary random variable $W^*$ with similar properties).
However, establishing such a convention does not seem realistic.

\medskip

\noindent Now it might be argued that the improvement in the properties of the confidence interval
for $\theta$ justifies the ``splitting'' of each observation $y$ and that such a split does
not require any meaningful comparison of the values that make up this split. However,
even if a convention could be enforced that the data-randomized
confidence intervals for $\theta$ are based
only on Korn's auxiliary random variable $W$, these confidence intervals would still not satisfy
 the
invariance property described in Example 2.35 on Cox and Hinkley (1974).

\bigskip

\noindent {\bf 7. Confidence intervals for $\theta$ based on splitting the Bernoulli data into two groups
of approximately equal relatively prime size}

\medskip

As before, suppose that $X_1, X_2,\ldots,X_n$ are independent and
identically Bernoulli$(\theta)$ distributed and that our aim is to
find a confidence interval for
$\theta$ that satisfies \eqref{upper_tail} and \eqref{lower_tail}.
In this section, we consider the properties of a confidence interval for $\theta$ that is obtained as follows.
Suppose that $n = n_1 + n_2$, where $n_1$ and $n_2$ are relatively prime and as close as possible.
Form the following estimator of $\theta$:
\begin{equation*}
\hat{\Theta} = \frac{1}{2} \left(\frac{Y_1}{n_1} + \frac{Y_2}{n_2} \right),
\end{equation*}
where $Y_1 = X_1 + \dots + X_{n_1}$ and $Y_2 = X_{n_1+1} + \dots + X_n$.
This is an unbiased estimator of $\theta$.
Consider the following procedure for finding a $1-\alpha$
confidence interval for $\theta$.
The Clopper-Pearson interval $\big[\ell_{CP}(Y), u_{CP}(Y) \big]$
is the intersection of upper and lower
$1-\alpha/2$ confidence intervals for $\theta$ that are based on
inverting the family of hypothesis tests using the test statistic
$Y$ (or, equivalently, the test statistic $Y/n$). We can construct an analogous confidence interval that is the
intersection of upper and lower
$1-\alpha/2$ confidence intervals for $\theta$ that are based on the test statistic
$\hat{\Theta}$. Let us denote this confidence interval by $\big[ \ell^{\dag}(\bX), u^{\dag}(\bX) \big ]$.
This confidence interval is obtained by deterministically splitting the data into two parts.
A random splitting of the data into two parts is considered by Decrouez and Hall (2013b).

For concreteness, consider the particular case that $n = 47$.  Form the following estimator of $\theta$:
\begin{equation*}
\hat{\Theta} = \frac{1}{2} \left(\frac{Y_1}{23} + \frac{Y_2}{24} \right),
\end{equation*}
where $Y_1 = X_1 + \dots + X_{23}$ and $Y_2 = X_{24} + \dots + X_{47}$.
We have obtained $Y_1$ and $Y_2$ by splitting the 47 Bernoulli trials
into two groups of approximately equal relatively prime
size.
How do the confidence intervals $\big[\ell_{CP}(Y), u_{CP}(Y) \big]$
and $\big[ \ell^{\dag}(\bX), u^{\dag}(\bX) \big ]$ compare?
We expect the estimator $\hat{\Theta}$
 to be a somewhat less efficient estimator
of $\theta$ than the maximum likelihood estimator $Y/n$. This is because
we give the same weight to the estimators $Y_1/23$ and $Y_2/24$, when the
more accurate estimator $Y_2/24$ should have been given a larger weight.
On the other hand, the estimator $\hat{\Theta}$ has $24 \times 25 - 1 = 599$ possible values,
whereas $Y/n$ has only 48 possible values. We therefore expect that the
$P_{\theta} \big(\theta < \ell^{\dag}(\bX) \big)$ and $P_{\theta} \big(\theta > u^{\dag}(\bX) \big)$
will tend to be closer to $\alpha/2$ than
$P_{\theta} \big(\theta < \ell_{CP}(Y) \big)$ and $P_{\theta} \big(\theta > u_{CP}(Y) \big)$,
respectively
(cf. Decrouez and Hall, 2013a). The fact
that the estimator $\hat{\Theta}$ has many more possible values than the
estimator $Y/n$ can also be expected to lead to a shortening of the confidence
intervals that will, to some extent, compensate or even overcome the
inefficiency of the estimator $\hat{\Theta}$ by comparison with the estimator
$Y/n$.

Observe, however, that the confidence interval $\big[ \ell^{\dag}(\bX), u^{\dag}(\bX) \big ]$
may be viewed as a data-randomized confidence interval for $\theta$. As suggested in Section 5,
we use the data-randomized confidence interval of Korn (1987) as the standard against which
we judge $\big[ \ell^{\dag}(\bX), u^{\dag}(\bX) \big ]$.
We expect $\big[ \ell^{\dag}(\bX), u^{\dag}(\bX) \big ]$ to have coverage and expected
length properties that are inferior to the data-randomized confidence interval of Korn (1987).
For a start, the upper endpoints of the confidence intervals of Korn (1987)
are based on a statistic that can take
$2^{47} \approx  1.407 \times 10^{14}$ values. This is much larger than the 599 possible values of the
statistic $\hat{\Theta}$, on which $\big[ \ell^{\dag}(\bX), u^{\dag}(\bX) \big ]$ is based.
In addition, the statistic $\hat{\Theta}$ orders the data in the wrong way.
The statistic on which the confidence interval is based should always take a larger value for observed
value $y=t+1$ than for $y=t$.
However, for $y_1=t$ and $y_2=0$ the observed value is $y=t$ and
$\hat{\theta} = t/46$, which exceeds $\hat{\theta} = (t+1)/48$ when $y_1=0$ and $y_2=t+1$
(so that the observed value is $y=t+1$) when $t > 23$.
Also, $\hat{\theta}$ takes the same value, $1/2$, for $(y_1,y_2)=(0,n_2)$ (so that $y=n_2$)
and $(y_1,y_2)=(n_1,0)$ (so that $y=n_1$). In the language of Kabaila and Lloyd (2006),
$\hat{\Theta}$ is {\sl not} a ``refinement'' of $Y$.
The fact that $\hat{\Theta}$ orders the data in the wrong way may be interpreted as just
another manifestation of the inefficiency of this estimator.
This means that if we are prepared to consider data-randomized confidence intervals then we
should be using the data-randomized confidence interval of Korn (1987) instead of
the confidence interval $\big[ \ell^{\dag}(\bX), u^{\dag}(\bX) \big ]$ based on $\hat{\Theta}$.

\bigskip

\noindent {\bf 8. Discussion}

\medskip

Various kinds of randomized, pseudorandomized and data-randomized ``equi-tailed'' confidence intervals for
the binomial probability, based on iid Bernoulli observations, have been reviewed.
Of course, randomization, pseudorandomization and data-randomization can be combined in
various ways. For example, we could combine randomization with data-randomization.
Undoubtedly,
such confidence intervals will continue to be of theoretical interest.

The standard confidence interval that satisfies the ``equi-tailed'' coverage constraints
described in the paper is the Clopper-Pearson interval, which is not randomized (or
pseudorandomized or data-randomized). Broadening the class of allowable interval estimators
to include either randomization, pseudorandomization or data-randomization (or a combination
of some of these) may be viewed
as allowing one to use an additional resource. The theoretical question is: How well is this
additional resource being used?
We have asked and answered this question in the case of the unusual confidence interval
for the binomial probability described in Section 7.
Of course, whether randomized, pseudorandomized or data-randomized confidence intervals will
ever be used in practice is open to question.



\bigskip

\noindent {\bf Acknowledgment}

\medskip

The author is grateful to Peter Hall for helpful discussions. Much of the work described in the paper
was carried out during the author's visit to the Department of Mathematics and Statistics, University
of Melbourne.


\bigskip

\noindent {\bf Appendix A: A generalization of the randomized confidence intervals}

\medskip

A generalization of the randomized confidence interval $\big[\ell_R(Y,V), u_R(Y,V) \big]$ is
$\big[\ell_R(Y,V_{\ell}), \newline u_R(Y,V_u) \big]$, where $Y$ and $(V_{\ell}, V_u)$ are independent,
$V_{\ell} \sim U(0,1)$ and $V_u \sim U(0,1)$.
One could, for example, choose $V_{\ell} = 1 - V_u$.
This interval satisfies the following conditions:
\begin{align*}
P_{\theta} \big(\theta > u_R(Y,V_u) \big ) &= \frac{\alpha}{2} \ \ \text{for all} \ \ \theta \\
P_{\theta} \big(\theta < \ell_R(Y,V_{\ell}) \big ) &= \frac{\alpha}{2} \ \ \text{for all} \ \ \theta.
\end{align*}
Also, the confidence intervals $\big[\ell_R(Y,V), u_R(Y,V) \big]$ and
$\big[\ell_R(Y,V_{\ell}), u_R(Y,V_u) \big]$
have the same expected length functions. There seems to be no advantage to be gained from this generalization.
For example, suppose that $V_{\ell} = 1 - V_u$.
In this case, the confidence interval lower endpoint $\ell_R(y,v_{\ell}) = \ell_R(y,1-v_u)$ is a decreasing
function of $v_u$ for (a) $y=0$ and $1-v_u \ge 1-\alpha/2$ and (b) all $y \ge 1$.
This means that, for given observed value $y$, the main effect of increasing $v_u$ is to widen the confidence interval.
In statistical practice, confidence interval width is interpreted as a measure of the accuracy of the estimation
of $\theta$. It does not seem helpful to report (according to this interpretation) varying apparent accuracies
of estimation of $\theta$ (depending on the value of $v_u$), for the same observed value $y$.

\bigskip

\noindent {\bf Appendix B: Randomized confidence intervals that depend on an auxiliary discrete random variable}

\medskip

Suppose that the random variable $W$ is such that, conditional on $Y=y$, $W$ is uniformly distributed
on $\big\{ 1/M(y), 2/M(y), \dots, M(y)/M(y) \big\}$, where $M(y)$ is an integer greater than 1 for each
$y = 0, \dots, n$. A particular case is that $M(y) = M$ for $y = 0, \dots, n$. Define the discrete random
variable $\tilde{W}$ by the requirement that, conditional on $Y=y$, $\tilde{W}=W-1/M(y)$. Thus, conditional on $Y=y$,
$\tilde{W}$ is uniformly distributed
on $\big\{0, 1/M(y), \dots, (M(y)-1)/M(y) \big\}$. Let $F_V$, $F_W$ and $F_{\tilde{W}}$ denote the cumulative
distribution functions of $V \sim U(0,1)$, $W$ and $\tilde{W}$, respectively. Using the facts that
$F_W$ is stochastically larger than $F_V$ and $F_V$ is stochastically larger than $F_{\tilde{W}}$,
it may be shown that the confidence interval $\big[\ell_R(Y,\tilde{W}), u_R(Y,W) \big]$
satisfies the following conditions:
\begin{align}
\label{upper_tail_discrete}
P_{\theta} \big(\theta > u_R(Y,W) \big ) &\le \frac{\alpha}{2} \ \ \text{for all} \ \ \theta \\
\label{lower_tail_discrete}
P_{\theta} \big(\theta < \ell_R(Y,\tilde{W}) \big ) &\le \frac{\alpha}{2} \ \ \text{for all} \ \ \theta.
\end{align}
If the smallest of the $M(y)$'s is not too small then, in addition,
\begin{align}
\label{upper_tail_discrete_approx}
P_{\theta} \big(\theta > u_R(Y,W) \big ) \approx \frac{\alpha}{2} \ \ \text{for all} \ \ \theta \\
\label{lower_tail_discrete_approx}
P_{\theta} \big(\theta < \ell_R(Y,\tilde{W}) \big ) \approx \frac{\alpha}{2} \ \ \text{for all} \ \ \theta
\end{align}
and the expected length functions of the confidence intervals
$\big[\ell_R(Y,V), u_R(Y,V) \big]$ and $\big[\ell_R(Y,\tilde{W}), u_R(Y,W) \big]$
are approximately equal. These results
may be interpreted as resulting from the fact that $W$
may be viewed as an approximation to $V$, which has a uniform
distribution on $(0,1)$.

\newpage


\noindent {\bf Appendix C: Confidence intervals that depend on an auxiliary deterministic sequence}

\medskip

Suppose that $Y_1, Y_2, \dots$ are
independent and
identically Binomial($n,\theta$) distributed. In other words, suppose that we carry out independent
repetitions of the statistical experiment that gives rise to $Y$.
Let $v_1, v_2, \dots$ be a deterministic sequence of real numbers such
that either $v_k \in [0,1)$ for every $k=1, 2, \dots$ or $v_k \in (0,1]$ for every $k=1, 2, \dots$.
Now suppose that
$v_1, v_2, \dots$ is uniformly distributed modulo 1, as defined by Kuipers and Niederreiter (1974).
For any given irrational number $\lambda$,
the sequence $v_k = \{n \lambda \}$, where $\{ a \}$ denotes the fractional part of $a$,
possesses these properties and
may be viewed as a pseudorandom sequence.
The van der Corput sequence (defined e.g. on page 127 of Kuipers and Niederreiter, 1974)
possesses these properties
and may be viewed as a quasi-random
sequence.
Suppose that,
given the observation $y_k$ of $Y_k$,
we compute the confidence interval $\big[\ell_R(y_k,v_k), u_R(y_k,v_k) \big]$.
We use the notation
\begin{equation*}
{\cal I}({\cal A}) =
\begin{cases}
1 &\text{if } {\cal A} \ \ \text{is true} \\
0 &\text{if } {\cal A} \ \ \text{is false}
\end{cases}
\end{equation*}
where ${\cal A}$ is an arbitrary statement. It may be shown that, for each $\theta$,
\begin{align*}
&\frac{1}{m} \sum_{k=1}^m {\cal I}\big(\theta > u_R(Y_k,v_k)\big)\ \ \text{converges almost surely to} \ \
\frac{\alpha}{2} \ \ \text{and} \\
&\frac{1}{m} \sum_{k=1}^m {\cal I}\big(\theta < \ell_R(Y_k,v_k)\big)\ \ \text{converges almost surely to} \ \
\frac{\alpha}{2}
\end{align*}
as $m \rightarrow \infty$. It may also be shown that, similarly, the long-run average lengths of the confidence intervals
$\big[\ell_R(Y,V), u_R(Y,V) \big]$ and $\big[\ell_R(y_k,v_k), u_R(y_k,v_k) \big]$ are the same.

Alternatively, we may suppose that
$w_1, w_2, \dots$ is a deterministic sequence of real numbers such
that
$w_1, w_2, \dots$ is a periodic sequence with period $N$, where
$(w_1, w_2, \dots, w_N)$ is a permutation of $(1/N, 2/N, \dots, N/N)$.
There are both pseudorandom sequences (found using e.g. mixed congruential generators
with maximal possible cycle length $N$) and sequences with very obvious pattern
(e.g. $1/N, 2/N, \dots, N/N, 1/N, 2/N, \dots$)
that satisfy these conditions.
Define the sequence $\tilde{w}_1, \tilde{w}_2, \dots$ by
$\tilde{w}_k = w_k - 1/N$ for $k = 1, 2, \dots$.
Suppose that,
given the observation $y_k$ of $Y_k$,
we compute the confidence interval $\big[\ell_R(y_k,\tilde{w}_k), u_R(y_k,w_k) \big]$.
It may be shown that, for each $\theta$,
\begin{align*}
&\frac{1}{m} \sum_{k=1}^m {\cal I}\big(\theta > u_R(Y_k,v_k)\big)\ \ \text{converges almost surely to a number} \ \
\le \frac{\alpha}{2} \ \ \text{and} \\
&\frac{1}{m} \sum_{k=1}^m {\cal I}\big(\theta < \ell_R(Y_k,v_k)\big)\ \ \text{converges almost surely to a number} \ \
\le \frac{\alpha}{2}
\end{align*}
as $m \rightarrow \infty$.

\bigskip

\noindent {\bf References}



\smallskip

\noindent Blaker, H., (2000). Confidence curves and improved exact confidence intervals for
discrete distributions. Canadian Journal of Statistics 4, 783--798.

\smallskip

\noindent Casella, G., Berger, R. L., (1990). {\sl Statistical
Inference}. Brooks/Cole, Pacific Grove, California.

\smallskip

\noindent Casella, G., Berger, R. L., (2002). {\sl Statistical
Inference, second edition}. Duxbury, Pacific Grove, California.

\smallskip

\noindent Clopper, C.J., Pearson, E.S., (1934). The used of confidence or fiducial limits illustrated
in the case of the binomial. Biometrika 26, 404--413.

\smallskip

\noindent Cox, D.R., Hinkley, D.V., (1974). {\sl Theoretical Statistics}. Chapman and Hall, London.

\smallskip

\noindent Decrouez, G., Hall, P., (2013a). Normal approximation and smoothness for sums of means of
lattice-valued random variables. Bernoulli, DOI: 10.3150/12-BEJSP02.

\smallskip

\noindent Decrouez, G., Hall, P., (2013b). Split-sample methods for constructing confidence intervals
for binomial and Poisson parameters. Technical Report, Department of Mathematics and Statistics,
Melbourne University.

\smallskip

\noindent Geyer, C.J., Meeden, G.D., (2005). Fuzzy and randomized confidence intervals and P-values.
Statistical Science 20, 358--366.

\smallskip

\noindent Kabaila, P., Lloyd, C.J., (2006). Improved Buehler confidence limits based on refined
designated statistics. Journal of Statistical Planning and Inference 136, 3145--3155.

\smallskip

\noindent Kiefer, J.K., (1987). {\sl Introduction to Statistical Inference}. Springer-Verlag, New York.

\smallskip

\noindent Korn, E.L., (1987). Data-randomized confidence intervals for discrete distributions.
Communications in Statistics: Theory and Methods 16,
705--715.

\smallskip

\noindent Kuipers, L., Niederreiter, H., (1974). {\sl Uniform Distribution of Sequences}.
John Wiley, New York.

\smallskip

\noindent Senn, S., (2007a). Drawbacks to noninteger scoring for ordered categorical data.
Biometrics 63, 269--299.

\smallskip

\noindent Senn, S., (2007b). {\sl Statistical Issues in Drug Development, second edition}. Wiley, Chichester, England.

\smallskip

\noindent Stevens, W. L., (1950). Fiducial limits of the parameter of a discontinuous random variable. Biometrika 37, 117--129.

\end{document}